\documentclass[12pt]{article}

\usepackage[left=3cm,right=3cm,top=4cm,bottom=4cm]{geometry}
\usepackage[T1]{fontenc}
\usepackage[latin1]{inputenc}
\usepackage{lmodern}
\usepackage[nottoc, notlof, notlot]{tocbibind}
\usepackage{amsmath}
\usepackage{amsthm}
\usepackage{amssymb}
\usepackage{ stmaryrd }
\usepackage{graphicx}
\usepackage{tikz} 
\usepackage{dsfont}
\usepackage{enumitem}
\usepackage{xspace}
\usepackage{ upgreek }
\usetikzlibrary{matrix,arrows} 
\usepackage[colorlinks,urlcolor=blue,linkcolor=blue,citecolor=brown]{hyperref}

\theoremstyle{plain}
\newtheorem{thm}{Theorem}

\theoremstyle{remark}

\theoremstyle{plain}
\newtheorem{prop}{Proposition}

\theoremstyle{definition}
\newtheorem*{defn}{Definition}

\theoremstyle{plain}
\newtheorem{lem}{Lemma}

\theoremstyle{plain}

\theoremstyle{plain}

\theoremstyle{remark}

\begin{document}

\title{There is no analogue to Jarn\'ik's relation for twisted Diophantine approximation}
\author{Antoine MARNAT\footnote{supported by IRMA Strasbourg, the Austrian Science Fund (FWF), Project F5510-N26, which is a part of the Special Research Program ``Quasi-Monte Carlo Methods: Theory and Applications'' and FWF START-project Y-901} \\
 \href{mailto:marnat@math.unistra.fr}{marnat@math.unistra.fr}\\ \\
}

\date{}
\maketitle

\abstract{ Jarn\'ik gave a relation between the two most classical uniform exponents of Diophantine approximation in dimension $2$. In this paper we consider a twisted case, between the classical and the multiplicative one, and we show that no analogue to Jarn\'ik's relation holds.  }



\section{Introduction and main result}

Given $\boldsymbol{\theta} = (\theta_1,\theta_2)$ with $1$, $\theta_1$, $\theta_2$ linearly independent over $\mathbb{Q}$, the exponent $\omega(\boldsymbol{\theta})$ (resp. the uniform exponent $\hat{\omega}(\boldsymbol{\theta})$) is defined as the supremum of the real numbers $\nu$ such that for arbitrarily large real number $H$ (resp.  for every sufficiently large real number $H$) the system of inequalities 
\[ 0<|q - p_1\theta_1 - p_2\theta_2| \leq H^{-\nu}, \; \; |p_1|\leq H \;  , \; \; |p_2|\leq H  \]
has an integer solution $(p_1,p_2,q)$.\\

On the other hand, the exponent $\lambda(\boldsymbol{\theta} )$ (resp. the uniform exponent $\hat{\lambda}(\boldsymbol{\theta} )$) is the supremum of the real numbers $\nu$ such that for arbitrarily large real number $H$ (resp. for every sufficiently large real number $H$) the system of inequalities
\[ 0< |q| \leq H, \; \; |q\theta_1   - p_1 |\leq H^{-\nu} \; , \; |q\theta_2  - p_2 |\leq H^{-\nu} \]
has an integer solution $(p_1,p_2,q)$.\\

Dirichlet's box principle (or Minkowski's first convex body theorem) provides the lower bounds

\[ \omega(\boldsymbol{\theta}) \geq \hat{\omega}(\boldsymbol{\theta}) \geq 2 \;  \textrm{ and } \; \lambda(\boldsymbol{\theta}) \geq \hat{\lambda}(\boldsymbol{\theta}) \geq \frac{1}{2}.   \]

The first result is Jarn\'ik's relation \cite{JAR} linking both uniform exponents:

\begin{thm}[Jarn\'ik, 1938]
For any pair of real numbers $\boldsymbol{\theta} = (\theta _1, \theta _2)$ with $1$, $\theta _1$, $\theta _2$ linearly independent over $\mathbb{Q}$, both uniform exponents satisfy the relation
\begin{equation}\tag{$\ast$}\label{RJ} \hat{\lambda}(\boldsymbol{\theta} ) + \frac{1}{\hat{\omega}(\boldsymbol{\theta} )} = 1. \end{equation}
\end{thm}

Recently,  Laurent \cite{ML} proved a theorem giving every possible value of the quadruple $\Omega(\boldsymbol{\theta} ) = (\omega(\boldsymbol{\theta} ), \lambda(\boldsymbol{\theta} ), \hat{\omega}(\boldsymbol{\theta} ), \hat{\lambda}(\boldsymbol{\theta} ))$ when $\boldsymbol{\theta} $ ranges over $\mathbb{R}^n$ such that $1$, $\theta_1$, $\theta_2$ are $\mathbb{Q}$-linearly independent, called \emph{spectrum}.

\begin{thm}[Laurent, 2009]\label{ML}
For any pair of real numbers $\boldsymbol{\theta} = (\theta _1, \theta _2)$ with $1$, $\theta _1$, $\theta _2$ linearly independent over $\mathbb{Q}$ the four exponents of Diophantine approximation satisfy the relations
\[\begin{array}{cccc}
2\leq \hat{\omega}(\boldsymbol{\theta}) \leq + \infty, & \hat{\lambda}(\boldsymbol{\theta}) + \frac{1}{\hat{\omega}(\boldsymbol{\theta})} =1, & \frac{\omega(\boldsymbol{\theta})(\hat{\omega}(\boldsymbol{\theta})-1)}{{\omega}(\boldsymbol{\theta})+\hat{\omega}(\boldsymbol{\theta})} \leq \lambda(\boldsymbol{\theta}) \leq \frac{\omega(\boldsymbol{\theta}) -(\hat{\omega}(\boldsymbol{\theta})-1)}{\hat{\omega}(\boldsymbol{\theta})}.\\
\end{array}\]
Conversely, for each quadruple $\omega, \hat{\omega}, \lambda,\hat{\lambda}$  satisfying 
\[\begin{array}{cccc}
2\leq \hat{\omega} \leq + \infty, & \hat{\lambda}+ \frac{1}{\hat{\omega}} =1, & \frac{\omega(\hat{\omega}-1)}{{\omega}+\hat{\omega}} \leq \lambda \leq \frac{\omega -(\hat{\omega}-1)}{\hat{\omega}},\\
\end{array}\]
 there exists a pair of real numbers $\boldsymbol{\theta} = (\theta _1, \theta _2)$ with $1$, $\theta _1$, $\theta _2$ linearly independent over $\mathbb{Q}$ with 
\[\begin{array}{cc}
\omega(\boldsymbol{\theta}) = \omega, & \lambda(\boldsymbol{\theta})= \lambda, \\
\hat{\omega}(\boldsymbol{\theta}) = \hat{\omega}, & \hat{\lambda}(\boldsymbol{\theta})= \hat{\lambda}.\end{array}\]
\end{thm}

In \cite{JAR}, Jarn\'ik noticed that there is no analogue to relation \eqref{RJ} in higher dimension. It is an open question to find an analogue to Jarn\'ik's relation in the multiplicative case, recently studied by German \cite{OG}.\\

Given $\boldsymbol{\theta} = (\theta_1,\theta_2)$ with $1$, $\theta_1$, $\theta_2$ linearly independent over $\mathbb{Q}$, the multiplicative exponent $\omega_{\times}(\boldsymbol{\theta})$ (resp. the uniform multiplicative exponent $\hat{\omega}_{\times}(\boldsymbol{\theta})$) is defined as the supremum of the real numbers $\nu$ such that for arbitrarily large real number $H$ (resp. for every sufficiently large real number $H$) the system of inequalities 
\[ 0 < |q - p_1\theta_1 - p_2\theta_2| \leq H^{-\nu}, \; \; \max(1,|p_1|) \max(1, |p_2|)\leq H^2 \]
has an integer solution $(p_1,p_2,q)$.\\

On the other hand, the multiplicative exponent $\lambda_{\times}(\boldsymbol{\theta})$ (resp. the uniform multiplicative exponent $\hat{\lambda}_{\times}(\boldsymbol{\theta})$) is the supremum of the real numbers $\nu$ such that for arbitrarily large real number $H$ (resp. for every sufficiently large real number $H$) the system of inequalities
\[ 0< |q| \leq H, \; \; |q\theta_1   - p_1 | \cdot |q\theta_2   - p_2 | \leq H^{-2\nu} \]
 has an integer solution $(p_1,p_2,q)$.\\
 
 Minkowski's first convex body theorem \cite{Mink} provides the lower bounds
 \[ \omega_\times(\boldsymbol{\theta}) \geq \hat{\omega}_\times(\boldsymbol{\theta}) \geq 2 \;  \textrm{ and } \; \lambda_\times(\boldsymbol{\theta}) \geq \hat{\lambda}_\times(\boldsymbol{\theta}) \geq \frac{1}{2}.   \]
 
We observe that
\[ \{ x,y \in \mathbb{R} \mid \max(1,|x|) \max(1, |y|) \leq H^2 \} = \cup_{0\leq i \leq 1} \{ x,y \in \mathbb{R} \mid |x| \leq H^{2i}, |y| \leq H^{2(1-i)} \},   \]
where we have a union of uncountably many sets. It is thus natural to consider the following twisted exponents. See the paper of Harrap \cite{Harrap} for further twisted Diophantine approximation.\\ 

Given $\boldsymbol{\theta} = (\theta_1,\theta_2)$ with $1$, $\theta_1$, $\theta_2$ linearly independent over $\mathbb{Q}$, and $i$, $j$ non negative real numbers with $i+j=1$,  the twisted exponent $\omega_{i,j}(\boldsymbol{\theta})$ (resp. the uniform twisted exponent $\hat{\omega}_{i,j}(\boldsymbol{\theta})$) is defined as the supremum of the real numbers $\nu$ such that for arbitrarily large real number $H$ (resp. for every sufficiently large real number $H$) the system of inequalities 
\[ 0<|q - p_1\theta_1 - p_2\theta_2| \leq H^{-\nu}, \; \; |p_1|\leq H^{2i}, \; \; |p_2|\leq H^{2j} \]
has an integer solution $(p_1,p_2,q)$ with $(p_1,p_2)\neq(0,0)$.\\

On the other hand, the twisted exponent $\lambda_{i,j}(\boldsymbol{\theta})$ (resp. the uniform twisted exponent $\hat{\lambda}_{i,j}(\boldsymbol{\theta})$) is the supremum of the real numbers $\nu$ such that for arbitrarily large real number $H$ (resp. for every sufficiently large real number $H$) the system of inequalities
\[ 0< |q| \leq H, \; \; |q\theta_1  - p_1 |\leq H^{-2i\nu}, \; \; |q\theta_2   - p_2 |\leq H^{-2j\nu}  \]
has an integer solution $(p_1,p_2,q)$.\\

 Again, Minkowski's first convex body theorem \cite{Mink} provides the lower bounds
 \[ \omega_{i,j}(\boldsymbol{\theta}) \geq \hat{\omega}_{i,j}(\boldsymbol{\theta}) \geq 2 \;  \textrm{ and } \; \lambda_{i,j}(\boldsymbol{\theta}) \geq \hat{\lambda}_{i,j}(\boldsymbol{\theta}) \geq \frac{1}{2}.   \]

In his PhD thesis, the author gives relations analogue to those of Theorem \ref{ML} in the twisted case. At present, there is no construction proving that these relations are best possible.\\

The goal of this paper is to show that there is no analogue to Jarn\'ik's relation in the twisted case. We first need to introduce the following notion.\\

Fix two real numbers $R>1$ and $\mu>1$, we call \emph{$(\mu,R)$-sequences} a pair of sequences of positive integers $(A_n)_{n\geq1}$ and $(B_n)_{n\geq1}$ satisfying the following properties:

\begin{enumerate}[label=(\roman*)]
\item{There exist disjoint finite non-empty sets of prime numbers $S$ and $T$ such that, for each $n\geq1$, the set $S$ contains all the prime factors of $A_n$ and the set $T$ contains all the prime factors of $B_n$.}
\item{For each $p\in S$ and $q\in T$, the sequence $(\nu_p(A_n))_{n\geq1}$ and $(\nu_q(B_n))_{n\geq1}$ are strictly increasing sequences of positive integers.}
\item{We have \[\lim_{n\to+\infty} \frac{\log(A_{n+1})}{\log(A_n)} = \lim_{n\to+\infty} \frac{\log(B_{n+1})}{\log(B_n)} = \mu    \; \textrm{  and  }  \; \lim_{n\to+\infty} \frac{\log(B_n)}{\log(A_n)}=R.\] }
\end{enumerate}

For fixed $R>1$ and $\mu>1$ and any choice of $S$ and $T$, one can construct such sequences. For example, suppose that $S=\{p_1, \ldots , p_k\}$ and $T=\{q_1, \ldots , q_l\}$, where $p_1, \ldots , p_k, q_1, \ldots , q_l$ are distinct prime numbers. Setting $P=p_1p_2\cdots p_k$ and $Q=q_1q_2\cdots q_l$,  we consider
\[ A_n = P^{\lfloor a \mu^n \log Q \rfloor} \; \textrm{  and  }  \;  B_n=Q^{  \lfloor a \mu^n  R \log P \rfloor},\]
where $\lfloor x \rfloor$ is the integer part of the real number $x$. Then the sequences $(A_n)_{n\geq1}$ and $(B_n)_{n\geq1}$ form a pair of $(\mu,R)$-sequences, provided that the real parameter $a$ is large enough to ensure property \emph{(ii)}.\\


\begin{thm}\label{TC}
Let $i$, $j$ be two non negative real numbers with $i\geq j$ and $i+j=1$. Fix two real numbers $\mu$ and $R$ satisfying the conditions:
\begin{equation} \label{cond1}
\mu >2, \; R < j(\mu-2), \; R > \frac{\mu}{\mu-2} \; \textrm{ and } R> \frac{\mu}{i(\mu-1)}. 
\end{equation}

For any $(\mu,R)$-sequences $(A_n)_{n\geq1}$ and $(B_n)_{n\geq1}$, we consider the pair of real numbers
\[\boldsymbol{\theta} =(\theta_1,\theta_2) = \left( \sum_{n\geq1} A_n^{-1}, \sum_{n\geq1} B_n^{-1} \right). \]
 Its twisted exponents satisfy
\begin{eqnarray*}
\hat{\omega}_{i,j}(\boldsymbol{\theta}) &=& \min\left(2j \frac{\mu-1}{R}, 2i\frac{\mu-1}{\mu}R  \right),\\
\hat{\lambda}_{i,j}(\boldsymbol{\theta})& =& \min\left(\frac{1}{2i}(1-\frac{R}{\mu-1}), \frac{1}{2j}(1- \frac{\mu}{(\mu-1)R} ) \right),\\
\omega_{i,j}(\boldsymbol{\theta}) &=& 2i(\mu-1).\\
\end{eqnarray*}
\end{thm}

Hence, for any $(\mu,R)$-sequences, one gets a point $\boldsymbol{\theta}$. Theorem \ref{TC} then computes the three exponents $\hat{\omega}_{i,j}(\boldsymbol{\theta})$, $\hat{\lambda}_{i,j}(\boldsymbol{\theta})$ and ${\omega}_{i,j}(\boldsymbol{\theta})$ for each $i$, $j$ for which \eqref{cond1} holds.\\

Note that in the case $i=j=1/2$, one gets 
\begin{eqnarray*}
\hat{\omega}(\boldsymbol{\theta}) &=& \min\left(\frac{\mu-1}{R}, \frac{\mu-1}{\mu}R  \right),\\
\hat{\lambda}(\boldsymbol{\theta})& =& \min\left(1-\frac{R}{\mu-1}, 1- \frac{\mu}{(\mu-1)R}  \right),\\
\end{eqnarray*}

which agrees with Jarn\'ik's relation \eqref{RJ}.\\

With a suitable choice of the parameters $\mu$ and $R$, we can deduce the following theorem as a corollary.

\begin{thm}\label{TF} Let $i$, $j$ be two non negative real numbers with $i> j$ and $i+j=1$. Let $\hat{w}$ be a real number with $\hat{w}>6i$ and \[\hat{\lambda} \in \left( \frac{1}{2i}\left( 1- \frac{2j}{\hat{w}} \right) , \min\left(\frac{1}{2i}, \frac{1}{2j}\left( 1- \frac{2i}{\hat{w}} \right)\right) \right).\] There exists uncountably many $\boldsymbol{\theta}$ such that $\hat{\omega}_{i,j}(\boldsymbol{\theta})=\hat{w}$ and $\hat{\lambda}_{i,j}(\boldsymbol{\theta})= \hat{\lambda}$.\\
\end{thm}

Thus, no analogue to Jarn\'ik's relation holds for the twisted Diophantine approximation. It is still an open question to find the optimal interval for $\hat{\lambda}_{i,j}(\boldsymbol{\theta})$.\\

\paragraph{Notation:} Given two sequences $(a_n)_{n\geq1}$ and $(b_n)_{n\geq1}$, we write $a_n \gg b_n$ (resp. $a_n \ll b_n$)  if, for sufficiently large $n$, there exists a real number $c_1$ such that $a_n \geq c_1 b_n$ (resp. there exists a real number $c_2$ such that $ a_n \leq c_2 b_n$). Finally, we write $a_n \asymp b_n$ if both $a_n \gg b_n$ and $a_n \ll b_n$.

\section{Sequence of minimal points and proof of Theorems  \ref{TC} and  \ref{TF}  }

The main tool to compute the exponents of Diophantine approximation and prove Theorem \ref{TC} is the notion of minimal points, as introduced by Davenport and Schmidt in \cite[§3]{DavSchm}  or Jarn\'ik in \cite{JAR}.\\

Let $L$ and $N$ be two functions from $\mathbb{Z}^l$ to $\mathbb{R}_+$ where $l$ is a positive integer. 

\begin{defn}
A sequence of \emph{$(L,N)$-minimal points} $(M_k)_{k\geq0} \in (\mathbb{Z}^l)^\mathbb{N}$ is a sequence such that \begin{enumerate}[label=(\roman*)]
\item{$(N(M_k))_{k\geq0}$ is an increasing sequence with $N(M_0) \geq1$ ,}
\item{$(L(M_k))_{k\geq0}$ is a decreasing sequence with $L(M_0) \leq 1$, }
\item{for every $k\geq0$ and every point $M\in \mathbb{Z}^l$, if $N(M) < N(M_{k+1})$ then $L(M) \geq L(M_k)$. }
\end{enumerate}
\end{defn}

We can normalize any functions $L$ and $N$ to satisfy the conditions $N(M_0) \geq1$ and $L(M_0) \leq 1$. These conditions are not restrictive. They are imposed to simplify the further use of minimal points. Note that there is no unicity, because there may exist $M_k$ and $M_k'$ such that $L(M_k)=L(M_k')$ and $N(M_k)=N(M_k')$ but $M_k \neq M_k'$. \\

Let $\boldsymbol{x}=(x_0,x_1,x_1)$ be an integer triple. Consider the following functions with parameters $i$, $j$ and $\boldsymbol{\theta}$
\begin{eqnarray*}
L_{\lambda}(\boldsymbol{x}) &=& 
\max\left(|x_0\theta_1-x_1|^{1/(2i)} ,  |x_0\theta_2-x_2|^{1/(2j)} \right), \\
N_{\lambda}(\boldsymbol{x}) &=& 
|x_0| ,\\
L_{\omega}(\boldsymbol{x}) &=& 
 | x_1\theta_1 + x_2\theta_2 -x_0|,  \\
N_{\omega}(\boldsymbol{x}) &=& 
\max\left( |x_1|^{1/(2i)}, |x_2|^{1/(2j)}    \right),
\end{eqnarray*}
where $|x|$ is the absolute value of the real number $x$. They are related to the twisted exponents through the following proposition.

\begin{prop}\label{fact}\label{MA}
Let $\boldsymbol{\theta}=(\theta_{1},\theta_{2})$ be a pair of real numbers and $0<j \leq i<1$ such that $i+j=1$. Let $(a_{n})_{n\geq1}$ (resp. $(b_{n})_{n\geq1}$) be a sequence of $(L_{\lambda},N_{\lambda})$-minimal points of $\boldsymbol{\theta}$  (resp. $(L_{\omega},N_{\omega})$-minimal points). We have the relations
\[\begin{array}{llll}
 \hat{\lambda}_{i,j}(\boldsymbol{\theta}) & =  \liminf_{n\to\infty}\left(-\frac{\log L_{\lambda}(a_n) }{\log N_{\lambda}(a_{n+1})}\right),  & \lambda_{i,j}(\boldsymbol{\theta}) & =\limsup_{n\to\infty}\left(-\frac{\log L_{\lambda}(a_n) }{\log N_{\lambda}(a_n)}\right), \\[3mm]
 \hat{\omega}_{i,j}(\boldsymbol{\theta}) & = \liminf_{n\to\infty}\left(-\frac{\log L_{\omega}(b_n)}{\log N_{\omega}(b_{n+1})}\right), & \omega_{i,j}(\boldsymbol{\theta}) & =  \limsup_{n\to\infty} \left(-\frac{\log L_{\omega}(b_n) }{\log N_{\omega}(b_{n})}\right) .\\
\end{array}\]
\end{prop}


We omit the proof of Proposition \ref{MA}. It is an easy consequence of the definitions of a sequence of minimal points and of the twisted exponents.\\

This proposition gives us the exact values of the two uniform exponents if we have enough information about the tail of a sequence of minimal points, since only the asymptotic behavior is significant. From now on, every statement is implicitly considered for sufficiently large $n$.\\

For every $(\mu,R)$-sequences $(A_{n})_{n\geq1}$ and $(B_{n})_{n\geq1}$ with $\mu$ and $R$ satisfying the conditions \eqref{cond1} from Theorem \ref{TC}, we set $\boldsymbol{\theta}=(\theta_1,\theta_2) = \left( \sum_{n\geq1} A_n^{-1}, \sum_{n\geq1}  B_n^{-1} \right)$. We also set for every integer $n\geq1$ the integers

\begin{equation}\label{seqA'}
A_{n}' = \sum_{k=1}^{n} A_{n} A_{k}^{-1}  \; \textrm{ and } \; B_{n}' = \sum_{k=1}^{n}B_{n} B_{k}^{-1} .\end{equation}

With this notation, we have

\begin{lem}\label{lemA}
Under the assumption of Theorem \ref{TC}, the sequence of primitive integer points
\begin{equation}\label{seqA}
(A_{1}',A_{1},0), (B'_{1},0,B_{1}), \ldots , (A'_{n},A_{n},0),(B'_{n},0,B_{n}),(A'_{n+1},A_{n+1},0),\ldots
\end{equation}
consists ultimately of a sequence of $(L_{\omega},N_{\omega})$-minimal points.
\end{lem}

\begin{lem}\label{lemB}
Under the assumption of Theorem \ref{TC}, consider 
\[  C_{n} =(A_{n}B_{n},A'_{n}B_{n},A_{n}B'_{n})  \textrm{ and } D_{n} =(A_{n+1}B_{n},A'_{n+1}B_{n},A_{n+1}B'_{n}).\]

Ultimately, $C_n$ and $D_n$ consist in $(L_{\lambda}, N_{\lambda})$-minimal points . We denote by $E_{n,1}, \ldots , E_{n,s_n}$ (resp. $F_{n,1}, \ldots , F_{n,t_n}$) the intermediate $(L_{\lambda}, N_{\lambda})$-minimal points between $C_n=E_{n,0}$ and $D_n=E_{n,s_n+1}$ (resp. $D_n=F_{n,0}$ and $C_{n+1}=F_{n,t_n+1}$). \\

For $0\leq k \leq s_n$ and $0\leq k' \leq t_n$, the intermediate $(L_{\lambda}, N_{\lambda})$-minimal points are of the shape 
\[E_{n,k}=\left(E_{n,k,0},E_{n,k,1},\frac{E_{n,k,0}}{B_n}B'_n\right),\]
 \[F_{n,k'}=\left(F_{n,k',0},\frac{F_{n,k',0}}{A_n}A'_n,F_{n,k',2}\right).\]  

Furthermore, they satisfy \\
\[\begin{array}{ll}
 L_{\lambda} (F_{n,k}) = \|F_{n,k,0}\theta_2\|^{\frac{1}{2j}} \asymp (\frac{1}{N_{\lambda}(F_{n,k+1})}A_{n+1})^{\frac{1}{2j}} ,&
 N_{\lambda}(F_{n,1}) \asymp \frac{B_{n+1}}{B_n} , \\
 L_{\lambda} (E_{n,k}) = \|E_{n,k,0}\theta_1\|^{\frac{1}{2i}} \asymp (\frac{1}{N_{\lambda}(E_{n,k+1})} B_n)^{\frac{1}{2i}}, &
  N_{\lambda}( E_{n,1}) \asymp \frac{A_{n+1}}{A_n}, \\
\end{array}\]
where $\|x\|$ denotes the distance from the real number $x$ to a nearest integer. 
\end{lem}

Lemmas \ref{lemA} and \ref{lemB} will be proved in Section \ref{prooflem}. We show now how Theorem \ref{TC} and \ref{TF} can be derived from Proposition \ref{MA} and Lemmas \ref{lemA} and \ref{lemB}.\\

{\it Proof of Theorem \ref{TC}}.\\

First, notice that 
\[|A_{n}\theta_{1}-A_{n}'| \asymp |A_{n}A_{n+1}^{-1}| \; \textrm{ and } |B_{n}\theta_{2}-B_{n}'| \asymp |B_{n}B_{n+1}^{-1}| .\]
Using Proposition \ref{MA}, we can thus compute

\begin{eqnarray*}
\hat{\omega}_{i,j}(\boldsymbol{\theta}) &=& \min\left( \liminf_{n\to\infty} \frac{-\log|A_{n}'-A_{n}\theta_{1}|}{\log|B_{n}|^{\frac{1}{2j}}} ,  \liminf_{n\to\infty} \frac{-\log|B_{n}'-B_{n}\theta_{2}|}{\log|A_{n+1}|^{\frac{1}{2i}}}   \right) \\
 &=& \min \left(  2j \liminf_{n\to\infty} \frac{\log A_{n+1} - \log A_{n}}{\log B_{n}} , 2i \liminf_{n\to\infty} \frac{\log B_{n+1}-\log B_{n}}{\log A_{n+1}}   \right) \\
&=& \min\left(2j \frac{\mu-1}{R}, 2i\frac{\mu-1}{\mu}R  \right),\end{eqnarray*}
\begin{eqnarray*}
\hat{\lambda}_{i,j}(\boldsymbol{\theta})& =&   \min\left(  \liminf_{n\to\infty} \min_{0\leq k \leq s_n}\left( -\frac{\log L_{\lambda}(E_{n,k}) }{\log N_{\lambda}(E_{n,k+1})} \right),  \liminf_{n\to\infty} \min_{0\leq k \leq t_n}\left(  -\frac{\log L_{\lambda}(F_{n,k}) }{\log N_{\lambda}(F_{n,k+1})} \right)   \right)    \\
&=&   \min\left(  \liminf_{n\to\infty} \min_{0\leq k \leq s_n} \frac{1}{2i} \left( 1 - \frac{\log B_n}{\log E_{n,k+1,0}} \right)  \liminf_{n\to\infty} \min_{0\leq k \leq t_n}\frac{1}{2j}\left( 1- \frac{\log A_{n+1}}{ \log F_{n,k+1,0} } \right)   \right) , \\
&=& \min\left(  \liminf_{n\to\infty} \frac{1}{2i} \left( 1 - \frac{\log B_n}{\log E_{n,1,0}} \right) ,  \liminf_{n\to\infty} \frac{1}{2j}\left( 1- \frac{\log A_{n+1}}{ \log F_{n,1,0} } \right)   \right)  \\
&=& \min\left(\frac{1}{2i}\left(1-\frac{R}{\mu-1}\right), \frac{1}{2j}\left(1- \frac{\mu}{(\mu-1)R} \right) \right),\end{eqnarray*}
\begin{eqnarray*}
{\omega}_{i,j}(\boldsymbol{\theta}) &=& \max\left( \limsup_{n\to\infty} \frac{-\log|A_{n}'-A_{n}\theta_{1}|}{\log|A_{n}|^{\frac{1}{2j}}} ,  \limsup_{n\to\infty} \frac{-\log|B_{n}'-B_{n}\theta_{2}|}{\log|B_{n}|^{\frac{1}{2i}}}   \right) \\
&=& \max\left( 2j \limsup_{n\to\infty}  \frac{\log A_{n+1} - \log A_{n}}{\log|A_{n}|}      , 2i \limsup_{n\to\infty} \frac{\log B_{n+1}-\log B_{n}}{\log|B_{n}|}   \right)  \\
&=& \max \left( 2j (\mu-1), 2i (\mu-1)  \right) = 2i (\mu-1) .  \\
\end{eqnarray*}

This completes the proof of Theorem \ref{TC}. \qed \\

{\it Proof of Theorem \ref{TF}}.\\

Theorem \ref{TF} is a corollary from Theorem \ref{TC} if we choose the parameters from the following proposition 

\begin{prop}\label{Propar}
For every $\hat{\omega}> 6i$, the parameters
\begin{eqnarray*}
R &\in & \left( \frac{\hat{\omega}}{2i}, \frac{\hat{\omega}}{2i}\frac{1+ \sqrt{1+\frac{16ij}{\hat{\omega}^{2}}}}{2}   \right), \\
\mu &=& \frac{2iR}{2iR-\hat{\omega}}.\\
\end{eqnarray*}
satisfy the conditions \eqref{cond1} from Theorem \ref{TC}.
\end{prop}

The key point is to notice that \[R \in  \left( \frac{\hat{\omega}}{2i}, \frac{\hat{\omega}}{2i}\frac{1+ \sqrt{1+\frac{16ij}{\hat{\omega}^{2}}}}{2}   \right)\] is equivalent to  \[ 0< R(2iR-\hat{\omega}) < 2j.\]

With the choices of the Proposition \ref{Propar}, we have

\begin{eqnarray*}
\hat{\omega}_{i,j}(\boldsymbol{\theta}) &=& \min\left( 2j \frac{\mu-1}{R}, 2i \frac{\mu-1}{\mu} R \right) \\
  &=& \min\left( 2j\frac{\hat{\omega}}{R(2iR-\hat{\omega})}, \hat{\omega}\right) = \hat{\omega},\\
  \hat{\lambda}_{i,j}(\boldsymbol{\theta}) &=& \min\left(   \frac{1}{2i}\left(1-\frac{R(2iR-\hat{\omega})}{\hat{\omega}}\right), \frac{1}{2j}\left( 1-\frac{2i}{\hat{\omega}}\right)\right)\\
  & \in & \left(  \frac{1}{2i}\left(1-\frac{2j}{\hat{\omega}}\right), \min\left( \frac{1}{2i}, \frac{1}{2j}\left(1-\frac{2i}{\hat{\omega}} \right)  \right) \right).
\end{eqnarray*}

To get uncountably many pairs $\boldsymbol{\theta}$ with given exponents, note that Theorem \ref{TC} also holds for \[\boldsymbol{\theta}=(\theta_1,\theta_2) = \left(\sum_{n\geq1} \varepsilon_n A_n^{-1}, \sum_{n\geq 1} \varepsilon'_n B_n^{-1}   \right),\]
where $\varepsilon_n$ and $\varepsilon'_n$ are $\pm 1$. In this case we define 

\begin{equation*}
A_{n}' = \sum_{k=1}^{n} \varepsilon_n A_{n} \varepsilon_k A_{k}^{-1}  \; \textrm{ and } \; B_{n}' = \sum_{k=1}^{n}\varepsilon_n' B_{n} \varepsilon_k' B_{k}^{-1} .\end{equation*}

\qed

\section{ Proof of the lemmas on minimal points}\label{prooflem}

First of all, observe that in the definition of $(\mu,R)$-sequences, the condition \emph{(ii)} on the disjointness of the sets $S$ and $T$ of prime numbers ensures that for every $n\geq1$,
\[\gcd(A_{n},B_{n}) = \gcd(A_{n},A_{n}') = \gcd(B_{n},B_{n}') = 1,\]
where $A_n'$ and $B_n'$ are given by \eqref{seqA'}. 
This ensures that the terms in the sequences in Lemmas \ref{lemA} and \ref{lemB} are primitive points.\\

{\it Proof of Lemma \ref{lemA} }.\\

 Since $R>1$ and $i\geq j$ we have $|A_n|^{\frac{1}{2i}} < |B_n|^{\frac{1}{2j}} $. Furthermore, the condition $ R<j(\mu-2)$ from \eqref{cond1} implies that $\frac{\mu}{R}>\frac{i}{j}$. This implies that $|B_n|^{\frac{1}{2j}} < |A_{n+1}|^{\frac{1}{2i}}$. Thus, the sequence \eqref{seqA} is increasing with respect to $N_{\omega}$. To prove that it consists ultimately of consecutive $(L_{\omega},N_{\omega})$-minimal points , it is enough to show that, if $n$ is large enough, there is no primitive integer point $\boldsymbol{x}=(x_0,x_1,x_2)$ either with 
\begin{equation}\label{cas1}
N_{\omega}(\boldsymbol{x}) = \max \{|x_1|^{\frac{1}{2i}},|x_2|^{\frac{1}{2j}} \} <  B_n^{\frac{1}{2j}} \textrm{ and } L_{\omega}(\boldsymbol{x}) = |x_0-x_1\theta_1-x_2 \theta_2| < |A_n'-A_n\theta_1|,
\end{equation}
or with
\begin{equation}\label{cas2}
N_{\omega}(\boldsymbol{x}) =\max\{|x_1|^{\frac{1}{2i}},|x_2|^{\frac{1}{2j}} \} < A_{n+1}^{\frac{1}{2i}} \textrm{ and } L_{\omega}(\boldsymbol{x}) = |x_0-x_1\theta_1-x_2 \theta_2| < |B_n'-B_n\theta_2|.
\end{equation}
Assuming that there is a primitive point $(x_0,x_1,x_2)$ satisfying \eqref{cas1}, one finds 
\begin{multline}\label{Cas1} A_nB_nx_0 - A_n'B_nx_1 - A_nB_n'x_2 = A_nB_n(x_0-x_1\theta_1-x_2\theta_2) \\- B_n(A_n'-A_n\theta_1)x_1 - A_n(B_n'-B_n\theta_2)x_2\end{multline}

The left hand side of this equality is an integer. If it is $0$, then $B_n$ divides $A_nB_n'x_2$, thus $B_n$ divides $x_2$ and so $x_2=0$ because of the hypothesis $|x_2|<B_n$. Thus we obtain $A_nx_0-A_n'x_1=0$. As $\boldsymbol{x}$ is a primitive point, this implies that $(x_0,x_1)= \pm(A_n',A_n)$ and so $|x_0 - x_1\theta_1-x_2\theta_2|=|A_n'-A_n\theta_1|$ contrary to the hypothesis. So the left hand side of \eqref{Cas1} has absolute value at least $1$ and thus
\begin{eqnarray*}
1 & \leq & A_nB_n|A_n'-A_n\theta_1| + B_n|A_n'-A_n\theta_1| |x_1| +A_n|B_n'-B_n\theta_2| |x_2|  \\ 
  & \ll &  A_n^2A_{n+1}^{-1}B_n + A_nA_{n+1}^{-1} B_n^{1+i/j} + A_nB_n^2B_{n+1}^{-1},
\end{eqnarray*}
because
\[|A_n'-A_n\theta_1| \leq 2\frac{A_n}{A_{n+1}} \;  \textrm{  and }  \; |B_n'-B_n\theta_2| \leq 2\frac{B_n}{B_{n+1}} .\]

We obtain a contradiction by showing that the three summands in the last expression above tend to zero as $n$ goes to infinity. First, $R< \mu-2$ implies that 
\[ R+2-\mu = \lim_{n\to\infty} \frac{2\log A_n + \log B_n - \log A_{n+1}}{\log  A_n} = -\rho_1 < 0. \]
That is,
\[  A_n^2A_{n+1}^{-1}B_n \sim_{n\to\infty} A_n^{-\rho_1} \to_{n\to_\infty} 0.\]
Then, $R<j(\mu-1)$ implies that
\[R+j(1-\mu) = \lim_{n\to\infty} \frac{j\log A_n + \log B_n - j \log A_{n+1}}{\log A_n} = -\rho_2 <0.\]
That is, 
\[A_nA_{n+1}^{-1} B_n^{1/j} \sim_{n\to \infty} A_n^{-\rho_2/j} \to_{n\to\infty} 0.\]
Finally, $(\mu-2)R >1$ implies that 
\[R(\mu-2)-1 = \lim_{n\to\infty} \frac{ \log A_n + 2 \log B_n - \log B_{n+1}}{-\log A_n} = \rho_3 >0.\]
That is, 
\[A_nB_n^2B_{n+1}^{-1} \sim_{n\to\infty} A_n^{-\rho_3} \to_{n\to\infty} 0.\]

Similarly, if the condition \eqref{cas2} has a solution $(x_0,x_1,x_2)$, $A_{n+1}$ can not divide $x_1$ and we deduce the lower bound
\begin{eqnarray*}
1 & \leq & A_{n+1}B_n|B_n'-B_n\theta_1| + B_n|A_{n+1}'-A_{n+1}\theta_1| |x_1| +A_{n+1}|B_n'-B_n\theta_2| |x_2|  \\ 
  & \ll &  B_n^2B_{n+1}^{-1}A_{n+1} + B_nA_{n+1}^2A_{n+2}^{-1} + A_{n+1}^{1+j/i}B_nB_{n+1}^{-1}.
\end{eqnarray*}

Again, the three summands in this expression tend to zero as $n$ goes to infinity because $R(\mu-2)>\mu$, $R< \mu(\mu-2)$ and $R(\mu-1) > \frac{\mu}{i}$ according to \eqref{cond1}.\\
 Thus, we proved Lemma \ref{lemA}. \qed \\
 
{\it Proof of  Lemma \ref{lemB}}.\\

 First, we show that ultimately the point $C_n$ and $D_n$ are $(L_{\lambda},N_{\lambda})$-minimal points, then we prove the properties of the $(L_{\lambda}, N_{\lambda})$-minimal points lying respectively between $C_n$ and $D_n$ and $D_n$ and $C_{n+1}$.\\

Since $A_{n+1} >A_n$ and $B_{n+1} >B_n$ by the property about valuations of $(\mu,R)$-sequences, the sequence
\[ C_1,D_1, \ldots  , C_n,D_n,C_{n+1}, \ldots \]
is strictly increasing with respect to $N_{\lambda}$.\\

Note that, if $n$ is large enough
\begin{eqnarray*}
L_{\lambda}(C_n) &=& \max \{  (|A_n\theta_1-A_n'| B_n)^{1/(2i)} ,  (|B_n\theta_2-B_n'| A_n)^{1/(2j)}   \} , \\
    &=& (|A_n\theta_1-A_n'| B_n)^{1/(2i)}, \\
    &\asymp& (A_nA_{n+1}^{-1} B_n)^{1/(2i)},
\end{eqnarray*}
since $|B_n\theta_2-B_n'| A_n \asymp A_nB_{n+1}^{-1} B_n$ goes faster to zero than $|A_n\theta_1-A_n'| B_n \asymp A_nA_{n+1}^{-1} B_n$ and we have $j \leq i$.\\

Suppose that a non-zero integer point $\boldsymbol{x}=(x_0,x_1,x_2)$ satisfies 
\begin{equation}\label{L}
L_{\lambda}(\boldsymbol{x}) \leq L_{\lambda}(C_n) \; \textrm{ and } \; N_{\lambda}(\boldsymbol{x}) =|x_0| \leq N_{\lambda}(D_n) = A_{n+1}B_n.
\end{equation}
We have the integer determinant
\begin{equation}\label{det1}
\left| \det\left( \begin{array}{cc} x_0 & x_2 \\ B_n & B_n' \end{array} \right) \right| \leq |x_0||B_n'-B_n\theta_2| + B_n|x_2-x_0\theta_2| \ll A_{n+1}B_n^2B_{n+1}^{-1} + B_nL_{\lambda}(C_n)^{2j}. \end{equation}

Since $R(\mu-2)  >\mu$,  the first summand $A_{n+1}B_n^2B_{n+1}^{-1}$ tends to zero as $n$ goes to infinity. Since $R<j(\mu-1)$ and $L_{\lambda}(C_n)\asymp (A_nA_{n+1}^{-1} B_n)^{1/(2i)}$, it follows that $B_nL_{\lambda}(C_n)^{2j}$ tends to zero as well. Thus, if $n$ is large enough, the determinant in \eqref{det1} is zero. Then, $(x_0,x_2)$ is a non-zero integral multiple of $(B_n,B_n')$, in particular $B_n$ divides $x_0$. Consider the integer determinant, for $\boldsymbol{x}=(x_0,x_1,x_2)$ with the stronger assumption that $N_{\lambda}(\boldsymbol{x}) \leq N_{\lambda}(C_n)=A_nB_n$:

\begin{equation}\label{det2}
\left| \det\left( \begin{array}{cc} x_0 & x_1 \\ A_n & A_n' \end{array} \right) \right| \leq |x_0||A_n'-A_n\theta_1| + A_n|x_1-x_0\theta_1| \ll A_{n}^2B_nA_{n+1}^{-1} + A_nL_{\lambda}(C_n)^{2i}. \end{equation}

Since $R < \mu-2$ summands tend to zero as $n$ goes to infinity. Thus for sufficiently large $n$, the determinant in \eqref{det2} is zero. Then, $(x_0,x_1)$ is a non-zero integral multiple of $(A_n,A_n')$, in particular $A_n$ divides $x_0$. Since $A_n$ and $B_n$ are coprime, $A_nB_n$ divides $x_0$ and $N_{\lambda}(\boldsymbol{x})\geq A_nB_n = N_{\lambda}(C_n)$. This proves that $C_n$ is a $(L_{\lambda},N_{\lambda})$-minimal point.\\

Similarly, if a non-zero integer point $\boldsymbol{x}=(x_0,x_1,x_2)$ satisfies
\[L_{\lambda}(\boldsymbol{x}) \leq L_{\lambda}(D_n) \asymp (A_{n+1}B_n B_{n+1}^{-1})^{1/(2j)} \; \textrm{ and } |x_0| \leq N_{\lambda}(C_{n+1}) = A_{n+1}B_{n+1}\]
and if $n$ is large enough, then $(x_0,x_1)$ is a non-zero integral multiple of $(A_{n+1},A_{n+1}')$ and if we strengthen the condition to $N_{\lambda}(\boldsymbol{x}) \leq N_{\lambda}(D_n)=A_{n+1}B_n$, then $(x_{0},x_{2})$ is a non zero integer multiple of $(B_n,B_n')$. Thus, $D_n$ is a $(L_{\lambda}, N_{\lambda})$-minimal point.\\

Consider any $(L_{\lambda}, N_{\lambda})$-minimal point $\boldsymbol{x}=(x_0,x_1,x_2)\in\mathbb{Z}^3$ with $N_{\lambda}(C_n) \leq x_0 < N_{\lambda}(D_n)$. We claim that for $n$ large enough, we have
\[ L_{\lambda}(\boldsymbol{x}) = |x_0\theta_1-x_1|^{1/(2i)}.  \]
Indeed, for such $\boldsymbol{x}$, recall that $(x_0,x_2)$ is an integer multiple of $(B_n,B_n')$. Since $x_0$ is strictly less than $A_{n+1}B_n$ and $B_n$ and $A_{n+1}$ are coprime, the points $(x_0,x_1)$ and $(A_{n+1},A_{n+1}')$ are linearly independent. Thus we have
\[ 1 \leq \left| \det\left( \begin{array}{cc} x_0 & x_1 \\ A_{n+1} & A_{n+1}' \end{array} \right) \right| \leq |x_0||A_{n+1}'-A_{n+1}\theta_1| + A_{n+1}|x_1-\theta_1x_0|. \]
This provides 
\begin{equation}\label{M1}
 |x_1-\theta_1x_0| \gg \frac{1}{A_{n+1}}, \end{equation}
 because $R > \mu(\mu-2)$ implies for sufficiently large $n$
 \[ |x_0||A_{n+1}'-A_{n+1}\theta_1| \leq  A_{n+1}^2B_nA_{n+2}^{-1} \leq \frac{1}{2}.\]
 
Furthermore, 
\begin{equation}\label{M2}
 |x_2-\theta_2x_0| = \frac{x_0}{B_n} |B_n\theta_2 - B_n'| \ll \frac{A_{n+1}B_n}{B_{n+1}}. \end{equation}
Finally, the condition $R>\frac{\mu}{i(\mu-1)}$ implies that $\left(\frac{A_{n+1}B_n}{B_{n+1}}\right)^{1/(2j)}$ tends to zero faster than $A_{n+1}^{-1/(2i)}$. Hence the result.\\

Let $E_{n,0}= C_{n}, E_{n,1}, \ldots , E_{n,s_n}, E_{n,s_n+1}=D_n$ denote a sequence of $(L_{\lambda},N_{\lambda})$-minimal points, chosen with positive first coordinate. Suppose that $n$ is large enough, by the above each $E_{n,k}$ has the form
\[ E_{n,k} = \left( E_{n,k,0}, E_{n,k,1}, \frac{E_{n,k,0}}{B_n}B_n' \right)  \;  \textrm{  where  }  \;  \frac{E_{n,k,0}}{B_n} \in \mathbb{N}^*  \;  (0\leq k \leq s_n+1), \]
and we have $L_{\lambda}(E_{n,k})^{2i}=|E_{n,k,0}\theta_1-E_{n,k,1}|$ for $k=0, \ldots , s_n$ but not for $k=s_n+1$.\\

We now show that \[ \frac{B_n}{2N_{\lambda}(E_{n,k+1})} \leq L_{\lambda}(E_{n,k})^{2i} \leq \frac{B_n}{ N_{\lambda}(E_{n,k+1})}. \]

By Minkowski's first convex body theorem, for each real number $H>0$, there exists a non-zero integer point $(x_0,x_1)$ satisfying the three conditions
\[ B_n \mid x_0,  \;  |x_0| \leq H  \;  \textrm{  and }  \; |x_0\theta_1-x_1|\leq \frac{B_n}{H}  .\]
For $n$ large enough, fix $k\in \{1, \ldots , s_n\}$ and choose $H$ in the range \[N_{\lambda}(E_{n,k}) < H <N_{\lambda}(E_{n,k+1}),\] then we get a non-zero integer point $\boldsymbol{x}=(x_0,x_1,x_2)$ where $x_0$ and $x_1$ are given by Minkowski's theorem and $x_2=x_0B_n^{-1}B_n'$. This point satisfies \[N_{\lambda}(E_{n,k}) < N_{\lambda}(\boldsymbol{x}) < N_{\lambda}(E_{n,k+1})\] and the conditions \eqref{M1} and \eqref{M2}. So,
\[ L_{\lambda}(\boldsymbol{x})^{2i} = |x_0\theta_1-x_1| \leq \frac{B_n}{H} .\]
By definition of $(L_{\lambda}, N_{\lambda})$-minimal points, we have $L_{\lambda}(E_{n,k}) < L_{\lambda}(\boldsymbol{x})$. Thus, by letting $H$ tend to $E_{n,k+1,0}$ we get the upper bound
\[ L_{\lambda}(E_{n,k})^{2i} \leq \frac{B_n}{N_{\lambda}(E_{n,k+1})}.\]

For the lower bound, notice that $B_n$ divides $E_{n,k,0}$ and $E_{n,k+1,0}$, and that two consecutive $(L_{\lambda}, N_{\lambda})$-minimal points are independent. Thus

\begin{eqnarray*}
B_n \leq \left| \det\left( \begin{array}{cc} E_{n,k,0} & E_{n,k,1} \\ E_{n,k+1,0} & E_{n,k+1,1} \end{array} \right) \right| &\leq & E_{n,k,0} |E_{n,k+1,1} - \theta_1E_{n,k+1,0} | + E_{n,k+1,0}|E_{n,k,1} - E_{n,k,0}\theta_1 | , \\
&\leq& 2 E_{n,k+1,0} |E_{n,k,1} - E_{n,k,0}\theta_1 |.
\end{eqnarray*}
So,
\[ L_{\lambda}(E_{n,k})^{2i} \geq \frac{B_n}{2N_{\lambda}(E_{n,k+1})}.\]

In particular, since $E_0=C_n$ we have
\[ N_{\lambda}(E_{n,1}) \asymp \frac{B_n}{L_{\lambda}(C_n)^{2i}} \asymp \frac{B_n}{A_nA_{n+1}^{-1} B_n } \asymp \frac{A_{n+1}}{A_n}. \]

Similarly, we show that for $n$ large enough, any $(L_{\lambda}, N_{\lambda})$-minimal point $\boldsymbol{x}=(x_0,x_1,x_2)$ with $N_{\lambda}(D_n) \leq N_{\lambda}(\boldsymbol{x}) < N_{\lambda}(C_{n+1})$  has 
\[ L_{\lambda}(\boldsymbol{x}) = |x_0\theta_2-x_2|^{1/(2j)},\]
and that the intermediate $(L_{\lambda}, N_{\lambda})$-minimal points $F_{n,k} = (F_{n,k,0}, \frac{F_{n,k,0}}{A_{n+1}}A_{n+1}', F_{n,k,2})$ satisfy 
\[ \frac{A_{n+1}}{2N_{\lambda}(F_{n,k+1})} \leq  L_{\lambda}(F_{n,k})^{2j} \leq \frac{A_{n+1}}{N_{\lambda}(F_{n,k+1})}.\]
via Minkowski's first convex body theorem and consideration on a good determinant.\\

This completes the proof of Lemma \ref{lemB}. \qed

\bibliographystyle{plain}
\bibliography{Biblio}

\begin{thebibliography}{1}

\bibitem{DavSchm}
Harold Davenport and Wolfgang~M. Schmidt.
\newblock Approximation to real numbers by algebraic integers.
\newblock {\em Acta Arithmetica}, 15:393--416, 1969.

\bibitem{OG}
Oleg~N. German.
\newblock Transference inequalities for multiplicative {D}iophantine exponents.
\newblock {\em Tr. Mat. Inst. Steklova}, 275(Klassicheskaya i Sovremennaya
  Matematika v Pole Deyatelnosti Borisa Nikolaevicha Delone):227--239, 2011.

\bibitem{Harrap}
Stephen Harrap.
\newblock Twisted inhomogeneous {D}iophantine approximation and badly
  approximable sets.
\newblock {\em Acta Arith.}, 151(1):55--82, 2012.

\bibitem{JAR}
Vojt{\v e}ch Jarn{\'\i}k.
\newblock Zum khintchineschen "{{\"U}}bertragungssatz".
\newblock {\em Trav. Inst. Math. Tbilissi}, 3:193--212, 1938.

\bibitem{ML}
Michel Laurent.
\newblock Exponents of diophantine approximmation in dimension two.
\newblock {\em Canad. J. Math.}, 61:165--189, 2009.

\bibitem{Mink}
Hermann Minkowski.
\newblock {\em Geometrie der {Z}ahlen}.
\newblock Bibliotheca Mathematica Teubneriana, Band 40. Johnson Reprint Corp.,
  New York-London, 1968.

\end{thebibliography}

\end{document}